\documentclass[11pt,a4paper]{article}
\usepackage[T1]{fontenc}
\usepackage{lmodern}
\usepackage{amsmath,amssymb,amsthm,mathtools}
\usepackage{booktabs,array}
\usepackage[margin=28mm]{geometry}
\usepackage{microtype}
\usepackage[colorlinks=true,
            citecolor=blue,
            linkcolor=blue,
            urlcolor=black,
            filecolor=black]{hyperref}
\usepackage{enumitem}
\usepackage{tikz}

\newtheorem{theorem}{Theorem}[section]
\newtheorem{lemma}[theorem]{Lemma}
\newtheorem{conjecture}{Conjecture}[section]
\newtheorem{problem}{Problem}[section]
\newtheorem{claim}[theorem]{Claim}

\newcommand{\ppath}{p_{\mathrm{path}}}
\newcommand{\bpath}{b_{\mathrm{path}}}
\newcommand{\pcyc}{p_{\mathrm{cyc}}}
\newcommand{\bcyc}{b_{\mathrm{cyc}}}

\setlist[enumerate]{label=\textup{(\roman*)},leftmargin=2.1em}

\title{Four-colour counterexamples for monochromatic path partitions}
\author{Zhen Liu\thanks{Email: \texttt{1552580575@qq.com}}
\and Chuanshu Wu\thanks{Email: \texttt{cswu97@126.com} (Corresponding author)}}
\date{\small Center for Discrete Mathematics, Fuzhou University,\protect\\
Fujian, 350003, China\protect\\[6pt]}

\begin{document}
\maketitle

\begin{abstract}
In 1989, Gyárfás conjectured that the vertex set of every \(r\)-edge-coloured complete graph can be partitioned into at most \(r\) monochromatic paths. The case \(r=2\) follows from a theorem of Gerencsér and Gyárfás, and Pokrovskiy proved the conjecture for \(r=3\). Recently, Li and Wang disproved the conjecture for all sufficiently large \(r\), leaving open the smallest number of colours for which it fails. We show that this number is four. For every integer \(n\geq118\), we construct a four-edge-colouring of \(K_n\) in which every monochromatic path partition contains at least five paths. We also construct, for every integer \(n\geq60\), a four-edge-colouring of \(K_{n,n}\) in which every monochromatic path partition contains at least eight paths, disproving Pokrovskiy’s conjectured bound \(2r-1\) for balanced complete bipartite graphs when \(r=4\).
\end{abstract}


\section{Introduction}

A monochromatic path partition of an edge-coloured graph $G$
is a collection of vertex-disjoint monochromatic paths whose
vertex sets partition $V(G)$. Monochromatic cycle partitions
are defined analogously. Paths are simple and nonempty, with
single vertices allowed. Single vertices and edges are also
regarded as cycles. Different paths or cycles in a partition
may have the same colour.
For $r\geq2$, let $\ppath(r)$ be the least integer $k$ such that
every $r$-edge-coloured complete graph admits a monochromatic
path partition into at most $k$ paths, and define $\pcyc(r)$
analogously for cycles. The study of these partitions goes back
to Gerencs\'er and Gy\'arf\'as~\cite{GG}, whose theorem gives
$\ppath(2)=2$; see the survey of Gy\'arf\'as~\cite{Gsurvey}.
Motivated by Rado's theorem for countably infinite complete
graphs~\cite{Rado}, Gy\'arf\'as posed the following conjecture
in 1989.

\begin{conjecture}[{\cite{G}}]\label{coj-Kn-path}
For every $r\geq3$, $\ppath(r)\leq r$.
\end{conjecture}

Pokrovskiy~\cite{P} proved that $\ppath(3)=3$, settling the
three-colour case.

Erd\H{o}s, Gy\'arf\'as and Pyber~\cite{EGP} proved
$\pcyc(r)=O(r^2\log r)$ and proposed the stronger conjecture
$\pcyc(r)\leq r$.
Gy\'arf\'as, Ruszink\'o, S\'ark\"ozy and Szemer\'edi~\cite{GRSS}
showed that $100r\log r$ cycles suffice for each fixed $r$
when the order of the complete graph is sufficiently large.
The two-colour cycle conjecture, due to Lehel, was established
through work of \L uczak, R\"odl and Szemer\'edi~\cite{LRS}
and Allen~\cite{Allen}, and finally proved for every order
by Bessy and Thomass\'e~\cite{BT}.
For three colours, three vertex-disjoint monochromatic cycles
cover all but $o(n)$ vertices~\cite{GRSS3}; Pokrovskiy~\cite{P23}
later reduced the number of uncovered vertices to an absolute
constant for sufficiently large $n$.
However, Pokrovskiy~\cite{P} constructed, for every $r\geq3$,
$r$-edge-coloured complete graphs whose vertex sets cannot be
partitioned into fewer than $r+1$ monochromatic cycles.
Related work includes path partitions in nearly complete
graphs~\cite{GJS}, cycle partitions in terms of the independence
number~\cite{Sarkozy}, and refinements of the connected-matching
method~\cite{Letzter}.

The analogous problem for balanced complete bipartite graphs
was studied by Gy\'arf\'as and Lehel~\cite{GL}
and Gy\'arf\'as~\cite{G83}.
As in~\cite{LW}, let $\bpath(r)$ and $\bcyc(r)$ denote the
corresponding path- and cycle-partition numbers for
$r$-edge-coloured balanced complete bipartite graphs.
Pokrovskiy~\cite{P} proved that $\bpath(2)=3$ and that
$\bpath(r)\geq2r-1$ for every $r\geq2$, and conjectured
the matching upper bound.

\begin{conjecture}[{\cite{P}}]\label{coj-bi-Kn-path}
For every $r\geq2$, $\bpath(r)\leq2r-1$.
\end{conjecture}

Together with Pokrovskiy's lower bound, a result of
Stein~\cite{Stein} gives
$5\leq\bpath(3)\leq6$.
For cycles, Haxell~\cite{Haxell} proved
$\bcyc(r)=O((r\log r)^2)$, which was improved to
$O(r^2\log r)$ by Peng, R\"odl and Ruci\'nski~\cite{PRR}.
Lang and Stein~\cite{LS} showed that $4r^2$ cycles suffice
for each fixed $r$ when the balanced complete bipartite graph
is sufficiently large.
For two colours, Schaudt and Stein~\cite{SS} gave a bound
of twelve cycles for sufficiently large graphs, and
Benevides, Quintino and Talon~\cite{BQT} proved
$\bcyc(2)\leq4$.
For three colours, Lang, Schaudt and Stein~\cite{LSS} proved
that five vertex-disjoint monochromatic cycles cover all but
$o(n)$ vertices, and that eighteen cycles suffice to partition
the entire vertex set for sufficiently large $n$.

Recently, Li and Wang~\cite{LW} disproved
Conjectures~\ref{coj-Kn-path} and~\ref{coj-bi-Kn-path}
for all sufficiently large $r$, proving
\[
 \ppath(r)\geq(1-o(1))r\log\log r
 \qquad\text{and}\qquad
 \bpath(r)=\Omega(r\log r).
\]
Their results leave open the smallest number of colours
for which Conjecture~\ref{coj-Kn-path} fails.
They asked the following question.

\begin{problem}[{\cite[Problem~5.3]{LW}}]
Determine the smallest $r$ for which $\ppath(r)>r$.
In particular, is $\ppath(4)>4$?
\end{problem}

We answer this question with the following theorem.

\begin{theorem}\label{thm:complete}
For every integer $n\geq118$, there exists a four-edge-colouring
of $K_n$ in which every monochromatic path partition contains
at least five paths.
\end{theorem}

Theorem~\ref{thm:complete} implies that $\ppath(4)\geq5$.
Together with $\ppath(2)=2$ and $\ppath(3)=3$, this shows that
four is the smallest number of colours for which
Conjecture~\ref{coj-Kn-path} fails.

We also give a four-colour construction for balanced complete
bipartite graphs.

\begin{theorem}\label{thm:bipartite}
For every integer $n\geq60$, there exists a four-edge-colouring
of $K_{n,n}$ in which every monochromatic path partition
contains at least eight paths.
\end{theorem}

Theorem~\ref{thm:bipartite} implies that $\bpath(4)\geq8$,
disproving Conjecture~\ref{coj-bi-Kn-path} for $r=4$.

\section{The complete graph case}\label{sec:complete}

\subsection{Construction}

\begin{figure}[htbp]
\centering
\definecolor{kNineOne}{HTML}{0072B2}
\definecolor{kNineTwo}{HTML}{D55E00}
\definecolor{kNineThree}{HTML}{009E73}
\definecolor{kNineFour}{HTML}{9B4F96}
\newcommand{\kNinePanel}[3]{%
  \begin{minipage}[t]{0.24\linewidth}
  \centering
  \begin{tikzpicture}[
      line width=0.75pt,
      line cap=round,
      line join=round
    ]
    \foreach \vertex [count=\index from 0]
      in {A,B,C,D,E,F,G,H,I}{
      \coordinate (\vertex) at ({90-40*\index}:1.25);
    }
    \foreach \first/\second/\third in {#2}{
      \draw[#1] (\first)--(\second)--(\third)--cycle;
    }
    \foreach \vertex [count=\index from 0]
      in {A,B,C,D,E,F,G,H,I}{
      \fill[white] (\vertex) circle[radius=0.055];
      \fill[black] (\vertex) circle[radius=0.035];
      \node[font=\small,inner sep=0pt]
        at ({90-40*\index}:1.49) {$\vertex$};
    }
  \end{tikzpicture}
  \par\smallskip
  {\small Colour #3}
  \end{minipage}%
}
\kNinePanel{kNineOne}{A/B/C,D/E/F,G/H/I}{1}\hfill
\kNinePanel{kNineTwo}{A/D/G,B/E/H,C/F/I}{2}\hfill
\kNinePanel{kNineThree}{A/E/I,B/F/G,C/D/H}{3}\hfill
\kNinePanel{kNineFour}{A/F/H,B/D/I,C/E/G}{4}
\caption{The four monochromatic subgraphs of $K_9$, each
consisting of three vertex-disjoint triangles.}
\label{fig:k9-colours}
\end{figure}
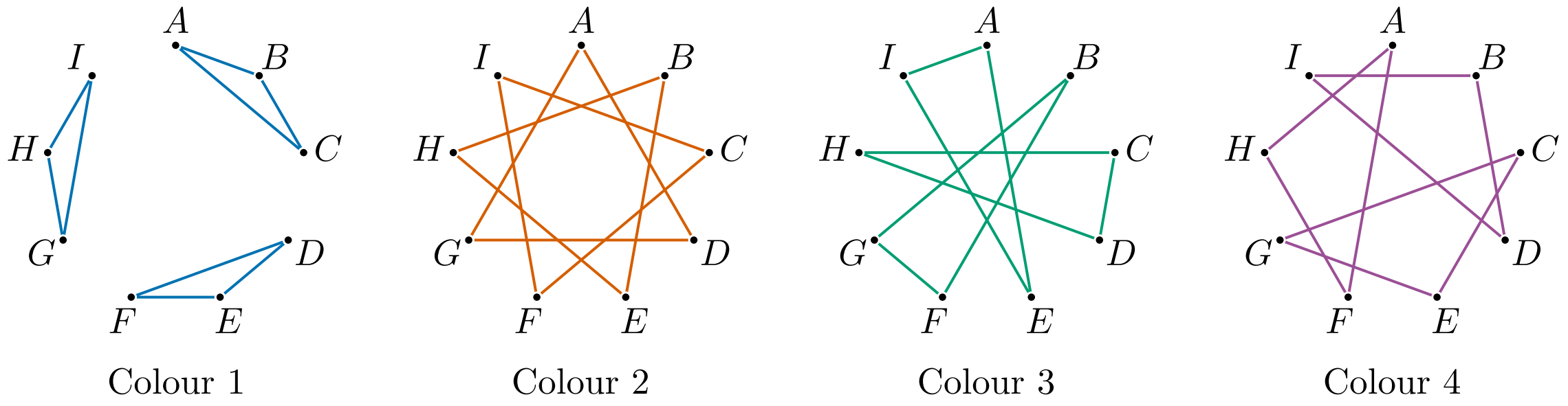

Let $V(K_9)=\{A,B,\ldots,I\}$. Colour the edges of the
following triangles with the indicated colours:
\[
\begin{array}{c@{\qquad}ccc}
\text{Colour }1: & ABC, & DEF, & GHI,\\
\text{Colour }2: & ADG, & BEH, & CFI,\\
\text{Colour }3: & AEI, & BFG, & CDH,\\
\text{Colour }4: & AFH, & BDI, & CEG.
\end{array}
\]
Here $ABC$ denotes the triangle with vertex set $\{A,B,C\}$.
Every edge of $K_9$ belongs to exactly one of the twelve
triangles above, so this defines a four-edge-colouring of
$K_9$. For each colour, the corresponding edges form three
vertex-disjoint triangles; see Figure~\ref{fig:k9-colours}.

Fix an integer $n\geq118$, and write
\[
n-118=3w+s,
\qquad
w\in\mathbb{Z}_{\geq0},
\quad
s\in\{0,1,2\}.
\]
Choose a partition
\[
V(K_n)=V_A\sqcup V_B\sqcup\cdots\sqcup V_I
\]
with class sizes given in Table~\ref{tab:weights}.
For distinct $X,Y\in\{A,B,\ldots,I\}$, colour every edge
between $V_X$ and $V_Y$ with the colour of $XY$ in the above
colouring of $K_9$. Colour the edges within each class with
the internal colour specified in Table~\ref{tab:weights}.

\begin{table}[htbp]
\centering
\begin{tabular}{c*{9}{c}}
\toprule
Class
& $V_A$ & $V_B$ & $V_C$ & $V_D$ & $V_E$
& $V_F$ & $V_G$ & $V_H$ & $V_I$\\
\midrule
Size
& $4+s$ & $32+w$ & $17$ & $1$ & $1$
& $4$ & $1$ & $35+w$ & $23+w$\\
Internal colour
& $3$ & $2$ & $1$ & $3$ & $3$
& $3$ & $3$ & $3$ & $4$\\
\bottomrule
\end{tabular}
\caption{Sizes and internal colours of the nine vertex classes.}
\label{tab:weights}
\end{table}

Since
\[
(4+s)+(32+w)+17+1+1+4+1+(35+w)+(23+w)
=118+3w+s=n,
\]
the required partition exists. This completes the construction
of the four-edge-colouring of $K_n$.
\subsection{Proof of Theorem~\ref{thm:complete}}

\begin{lemma}\label{lem:fourtriangles}
In the above colouring of $K_9$, three monochromatic triangles
cover all nine vertices only if they are the three triangles
of one colour. Four monochromatic triangles cover all nine
vertices only if they contain the three triangles of one
colour or have four different colours and a common vertex.
Repetitions are allowed.
\end{lemma}
\begin{proof}
We claim that $|V(T_1)\cap V(T_2)|=1$ for any two
monochromatic triangles $T_1,T_2$ of different colours.
Let $T_1',T_1''$ be the other two triangles of the same
colour as $T_1$.
Since triangles of different colours share no edge and
$V(K_9)=V(T_1)\sqcup V(T_1')\sqcup V(T_1'')$, we have
\[
1\geq |V(T_1)\cap V(T_2)|
=3-|V(T_1')\cap V(T_2)|-|V(T_1'')\cap V(T_2)|
\geq1.
\]

Suppose there exist three monochromatic triangles
$T_1,T_2,T_3$ such that
$\bigcup_{i=1}^{3}V(T_i)=V(K_9)$.
Since $\sum_{i=1}^{3}|V(T_i)|=9$, we have
$V(T_i)\cap V(T_j)=\varnothing$ for $i\ne j$.
By the preceding claim, any two monochromatic triangles of different colours intersect in exactly one vertex. It follows that $T_1, T_2, T_3$ share the same colour.

Now suppose $\bigcup_{i=1}^{4}V(T_i)=V(K_9)$ and
$T_1,\ldots,T_4$ do not include all three triangles of any colour.
By the preceding case, $T_1,\ldots,T_4$ are distinct, since any
repetition would leave three triangles whose vertex sets cover
$V(K_9)$, forcing them to be the three triangles of one colour.
If two have the same colour, say $T_1,T_2$, then
$V(T_1)\cap V(T_2)=\varnothing$, and neither $T_3$ nor $T_4$
has this colour. The claim gives
$|V(T_i)\setminus(V(T_1)\cup V(T_2))|=1$ for $i=3,4$.
Hence
\[
9=\left|\bigcup_{i=1}^{4}V(T_i)\right|
\leq6+1+1=8,
\]
a contradiction.
Thus $T_1,\ldots,T_4$ have distinct colours.
Set $U=V(T_1)\cup V(T_2)$. Since $T_1$ and $T_2$ have
different colours, the preceding claim gives
$|V(T_1)\cap V(T_2)|=1$; write
$V(T_1)\cap V(T_2)=\{x\}$. Thus $|U|=3+3-1=5$.
For $i\in\{3,4\}$, the triangle $T_i$ has a different colour
from both $T_1$ and $T_2$, so
$|V(T_i)\cap V(T_1)|=|V(T_i)\cap V(T_2)|=1$.
In particular, $|V(T_i)\setminus U|\leq2$. Hence
\[
9=\left|\bigcup_{i=1}^{4}V(T_i)\right|
\leq |U|+\sum_{i=3}^{4}|V(T_i)\setminus U|
\leq5+2+2=9.
\]
Equality implies $|V(T_i)\cap U|=1$ for $i=3,4$.
This unique vertex lies in both $V(T_1)$ and $V(T_2)$, and is
therefore $x$. Consequently,
$\bigcap_{i=1}^{4}V(T_i)=\{x\}$.
\end{proof}

\begin{proof}[Proof of Theorem~\ref{thm:complete}]
Suppose, to the contrary, that four pairwise vertex-disjoint
monochromatic paths $P_1,P_2,P_3,P_4$ cover $V(K_n)$.
Give each singleton path an arbitrary colour.

We first record a simple path inequality. Let
$P=v_1v_2\cdots v_\ell$ be a path of colour $c$, and let
$W\subseteq V(K_n)$ be such that no two vertices of $W$ are
joined by an edge of colour $c$. Since $v_iv_{i+1}$ has colour
$c$, the vertices $v_i$ and $v_{i+1}$ cannot both lie in $W$
for any $1\leq i<\ell$. Hence
\begin{equation}\label{eq:path-bound}
|V(P)\cap W|\leq |V(P)\setminus W|+1.
\end{equation}

For later use, we list the triples $(R;S,T)$ below, where
$R$ is distinguished and the order of $S,T$ is irrelevant.
For each listed triple, $RST$ is a triangle of colour $c$,
no edge within $V_R$ has colour $c$, and
$|V_R|\geq |V_S|+|V_T|+2$:
\begin{equation}\label{eq:obstruction-list}
\begin{array}{c@{\quad}l}
c=1: & (B;A,C),\ (F;D,E),\ (H;G,I),\\
c=2: & (A;D,G),\ (H;B,E),\ (I;C,F),\\
c=3: & (B;F,G),\ (I;A,E),\\
c=4: & (B;D,I),\ (C;E,G),\ (H;A,F).
\end{array}
\end{equation}
These inequalities follow directly from
Table~\ref{tab:weights}, using $w\geq0$ and $s\leq2$.

\begin{claim}\label{clm:class-obstruction}
Let $(R;S,T)$ be a triple listed for $c$ in
\eqref{eq:obstruction-list}. If $P$ is a path of colour $c$
with $V(P)\subseteq V_R\cup V_S\cup V_T$, then
$V_R\nsubseteq V(P)$.
\end{claim}

\begin{proof}
Since no edge within $V_R$ has colour $c$,
\eqref{eq:path-bound}, applied with $W=V_R$, yields
\[
|V(P)\cap V_R|
\leq |V(P)\cap(V_S\cup V_T)|+1
\leq |V_S|+|V_T|+1<|V_R|.
\]
Thus $V_R\nsubseteq V(P)$.
\end{proof}

We shall also use the following identities, which follow
directly from Table~\ref{tab:weights}:
\begin{equation}\label{eq:size-identities}
\begin{aligned}
|V_B|=|V_F|+|V_G|+|V_D|+|V_I|+3,
|V_I|+|V_C|=|V_G|+|V_D|+|V_H|+3.
\end{aligned}
\end{equation}

\begin{claim}\label{clm:supporting-triangles}
For each $i\in\{1,2,3,4\}$, there is a monochromatic triangle
$R_iS_iT_i$ of the same colour as $P_i$ such that
$V(P_i)\subseteq V_{R_i}\cup V_{S_i}\cup V_{T_i}$.
The four triangles $R_iS_iT_i$, $i\in\{1,2,3,4\}$, cover
$V(K_9)$.
\end{claim}

\begin{proof}
Fix $i\in\{1,2,3,4\}$, and let $c$ be the colour of $P_i$.
Choose $R_i\in\{A,\ldots,I\}$ such that
$V(P_i)\cap V_{R_i}\neq\varnothing$, and let $R_iS_iT_i$ be
the unique triangle of colour $c$ containing $R_i$.

Since the colour-$c$ subgraph of $K_9$ is the disjoint union
of three triangles, no edge of colour $c$ joins
$V_{R_i}\cup V_{S_i}\cup V_{T_i}$ to any other class.
As $P_i$ has colour $c$ and meets $V_{R_i}$, it follows that
$V(P_i)\subseteq V_{R_i}\cup V_{S_i}\cup V_{T_i}$.
We call $R_iS_iT_i$ the supporting triangle of $P_i$.

Now let $X\in\{A,\ldots,I\}$. Since $V_X\neq\varnothing$ and
the four paths cover $V(K_n)$, some $P_i$ contains a vertex of
$V_X$. The preceding inclusion then gives
$X\in\{R_i,S_i,T_i\}$. Hence the four supporting triangles
cover $V(K_9)$.
\end{proof}

By Claim~\ref{clm:supporting-triangles}, the four supporting
triangles cover $V(K_9)$. Lemma~\ref{lem:fourtriangles}
therefore implies that either they contain the three triangles
of one colour, or they have four distinct colours and a common
vertex. We consider these two cases separately.

\medskip
\noindent\textbf{Case 1.}
Three of the four supporting triangles are precisely the three
triangles of some colour $c$.

After relabelling the paths, assume that
$R_iS_iT_i$, $i\in\{1,2,3\}$, are the three triangles of
colour $c$. Thus $R_4S_4T_4$ is the remaining supporting
triangle.

\begin{claim}\label{clm:fourth-triangle}
If $(R;S,T)$ is listed for $c$ in
\eqref{eq:obstruction-list}, then
$R\in\{R_4,S_4,T_4\}$.
\end{claim}

\begin{proof}
Suppose that $R\notin\{R_4,S_4,T_4\}$. Since $RST$ is the
unique colour-$c$ triangle containing $R$, it is the only
supporting triangle containing $R$. By
Claim~\ref{clm:supporting-triangles}, any path containing a
vertex of $V_R$ must have a supporting triangle containing
$R$. Thus only the path supported by $RST$ can contain
vertices of $V_R$. Since the four paths cover $V(K_n)$, this
path contains all vertices of $V_R$, contradicting
Claim~\ref{clm:class-obstruction}.
\end{proof}

For $c=1$, Claim~\ref{clm:fourth-triangle} gives
$\{B,F,H\}\subseteq\{R_4,S_4,T_4\}$, and hence
$R_4S_4T_4=BFH$. This is impossible, since $BFH$ is not a
monochromatic triangle. Similarly, $c=2$ and $c=4$ would give
$R_4S_4T_4=AHI$ and $R_4S_4T_4=BCH$, respectively, neither of
which is monochromatic. Therefore $c=3$.
Claim~\ref{clm:fourth-triangle} now gives
$B,I\in\{R_4,S_4,T_4\}$. Since $BDI$ is the unique
monochromatic triangle containing the edge $BI$, we have
$R_4S_4T_4=BDI$.

Let $P$ and $Q$ be the paths supported by $BFG$ and $BDI$,
respectively. These are the only supporting triangles
containing $B$. Hence, by
Claim~\ref{clm:supporting-triangles}, the paths $P$ and $Q$
together cover $V_B$. The paths have colours $3$ and $4$,
whereas every edge within $V_B$ has colour $2$. Applying
\eqref{eq:path-bound} to each path with $W=V_B$, we obtain
\[
\begin{aligned}
|V_B|=|V(P)\cap V_B|+|V(Q)\cap V_B|\leq |V_F|+|V_G|+1+|V_D|+|V_I|+1=|V_B|-1,
\end{aligned}
\]
where the last equality follows from
\eqref{eq:size-identities}. This contradiction completes
Case~1.

\medskip
\noindent\textbf{Case 2.}
The four supporting triangles have distinct colours and a
common vertex $X$.

Any two monochromatic triangles of different colours intersect
in exactly one vertex. Hence any two supporting triangles
intersect precisely in $X$. Since they cover $V(K_9)$, every
vertex other than $X$ belongs to exactly one supporting
triangle.

\begin{claim}\label{clm:common-vertex}
We have $X=H$.
\end{claim}

\begin{proof}
We first show that $X\in\{B,E,H\}$. There is nothing to prove
if $X=B$, so suppose that $X\neq B$. Then exactly one
supporting triangle contains $B$. By
Claim~\ref{clm:supporting-triangles}, only the path supported
by this triangle can meet $V_B$. Since the four paths cover
$V(K_n)$, this path must contain all vertices of $V_B$.
Claim~\ref{clm:class-obstruction} excludes the triangles
$ABC$, $BFG$, and $BDI$, corresponding respectively to
$(B;A,C)$, $(B;F,G)$, and $(B;D,I)$. Thus the supporting
triangle containing $B$ is $BEH$. Since it also contains the
common vertex $X$, we obtain $X\in\{E,H\}$, and hence
$X\in\{B,E,H\}$ in all cases.

We next show that $X\in\{C,D,H\}$. Again, there is nothing to
prove if $X=H$. Suppose that $X\neq H$. Then exactly one
supporting triangle contains $H$, and its path must contain all
vertices of $V_H$. Claim~\ref{clm:class-obstruction} excludes
$GHI$, $BEH$, and $AFH$, corresponding respectively to
$(H;G,I)$, $(H;B,E)$, and $(H;A,F)$. Thus this triangle is
$CDH$, and so $X\in\{C,D,H\}$.

Consequently,
$X\in\{B,E,H\}\cap\{C,D,H\}=\{H\}$, proving the claim.
\end{proof}

The four supporting triangles are therefore
$GHI$, $BEH$, $CDH$, and $AFH$. Let $P$ and $Q$ be the paths
supported by $GHI$ and $CDH$, respectively. Since no other
supporting triangle contains $I$, Claim~\ref{clm:supporting-triangles}
implies that $P$ contains all vertices of $V_I$. Similarly,
$Q$ contains all vertices of $V_C$.

The path $P$ has colour $1$, whereas every edge within $V_I$
has colour $4$. Applying \eqref{eq:path-bound} with $W=V_I$
gives
\[
|V(P)\cap V_H|\geq |V_I|-|V_G|-1.
\]
Similarly, $Q$ has colour $3$, whereas every edge within
$V_C$ has colour $1$, and hence
\[
|V(Q)\cap V_H|\geq |V_C|-|V_D|-1.
\]
Since $P$ and $Q$ are vertex-disjoint,
\[
\begin{aligned}
|V_H|\geq |V(P)\cap V_H|+|V(Q)\cap V_H|\geq |V_I|+|V_C|-|V_G|-|V_D|-2=|V_H|+1,
\end{aligned}
\]
where the last equality follows from
\eqref{eq:size-identities}. This contradiction completes
Case~2.

Both cases are impossible. Therefore no four pairwise
vertex-disjoint monochromatic paths cover $V(K_n)$.
\end{proof}

\section{Balanced complete bipartite graphs}
\label{sec:bipartite}

\subsection{Construction}

Fix an integer $n\geq60$, and let $X,Y$ be the partite
sets of $K_{n,n}$. Partition $Y$ into twelve classes
$W_{ab}$, where $a,b\in\{1,2,3,4\}$ and $a\neq b$.
Give each class five vertices, and distribute the remaining
$n-60$ vertices as equally as possible among
$W_{12},W_{23},W_{34},W_{41}$. Thus every $W_{ab}$
has at least five vertices.

Partition $X$ into four classes with
\[
|V_1|=2\left\lceil\frac n4\right\rceil+3,\qquad
|V_2|=n-2\left\lceil\frac n4\right\rceil-5,\qquad
|V_3|=|V_4|=1.
\]
These sizes satisfy
$|V_1|>2(\lceil n/4\rceil+1)$ and
$|V_2|>\lceil n/4\rceil+1$; the latter follows from
$n-3\lceil n/4\rceil-6\geq(n-33)/4>0$.

\begin{figure}[htbp]
\centering
\definecolor{kBipOne}{HTML}{0072B2}
\definecolor{kBipTwo}{HTML}{D55E00}
\definecolor{kBipThree}{HTML}{009E73}
\definecolor{kBipFour}{HTML}{9B4F96}

\newcommand{\kBipPanel}[3]{%
  \begin{minipage}[t]{0.24\linewidth}
  \centering
  \begin{tikzpicture}[
      font=\footnotesize,
      line width=0.75pt,
      line cap=round,
      line join=round
    ]
    \path[use as bounding box]
      (-0.55,-4.10) rectangle (2.20,1.05);

    \foreach \i/\y in {
      1/0,
      2/-1.20,
      3/-2.40,
      4/-3.60
    }{
      \coordinate (V\i) at (0,\y);
    }

    \node[font=\small,inner sep=0pt]
      at (0,0.85) {$X$};
    \node[font=\small,inner sep=0pt]
      at (1.45,0.85) {$Y$};

    \foreach \i/\first/\second/\third in {#2}{
      \foreach \pair [count=\j from 0]
        in {\first,\second,\third}{
        \path (V\i)
          ++(1.45,{0.32-0.32*\j})
          coordinate (W\pair);
        \draw[#1] (V\i)--(W\pair);
        \fill[white] (W\pair) circle[radius=0.055];
        \fill[black] (W\pair) circle[radius=0.035];
        \node[anchor=west,inner sep=0pt,xshift=4pt]
          at (W\pair) {$W_{\pair}$};
      }
    }

    \foreach \i in {1,2,3,4}{
      \fill[white] (V\i) circle[radius=0.055];
      \fill[black] (V\i) circle[radius=0.035];
      \node[anchor=east,inner sep=0pt,xshift=-4pt]
        at (V\i) {$V_{\i}$};
    }
  \end{tikzpicture}
  \par\smallskip
  {\small Colour #3}
  \end{minipage}%
}

\kBipPanel{kBipOne}{
  1/12/13/14,
  2/21/31/41,
  3/23/24/34,
  4/32/42/43
}{1}\hfill
\kBipPanel{kBipTwo}{
  1/21/23/24,
  2/12/32/42,
  3/13/14/43,
  4/31/34/41
}{2}\hfill
\kBipPanel{kBipThree}{
  1/31/32/34,
  2/13/23/43,
  3/12/41/42,
  4/14/21/24
}{3}\hfill
\kBipPanel{kBipFour}{
  1/41/42/43,
  2/14/24/34,
  3/21/31/32,
  4/12/13/23
}{4}

\caption{The four monochromatic subgraphs in the bipartite
construction.}
\label{fig:bipartite-colours}
\end{figure}
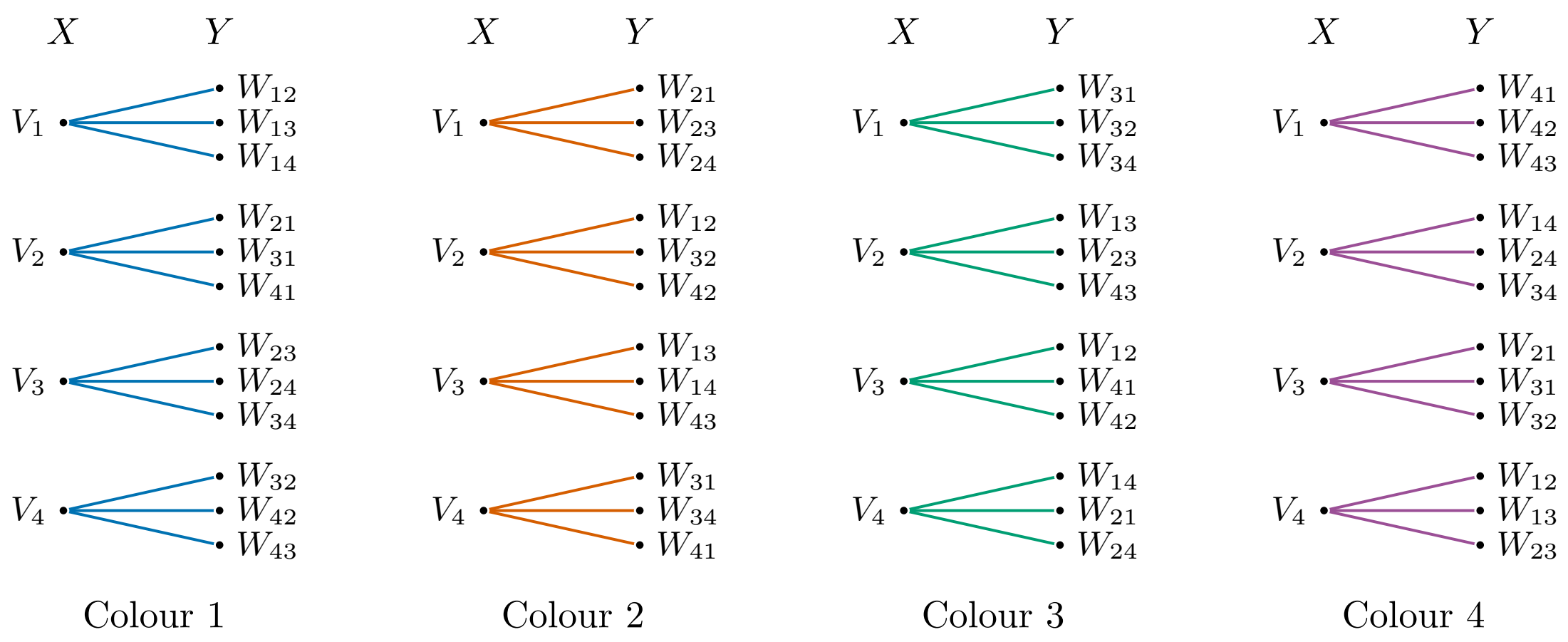

For each $y\in W_{ab}$, colour its edges to $V_1,V_2$
with $a,b$, respectively. Use the remaining two colours
on its edges to $V_3,V_4$, in increasing order if $a<b$
and in decreasing order if $a>b$.
The four monochromatic subgraphs are shown schematically
in Figure~\ref{fig:bipartite-colours}.

\subsection{Proof of Theorem~\ref{thm:bipartite}}

\begin{proof}[Proof of Theorem~\ref{thm:bipartite}]

For each $i\in\{1,2\}$ and each colour $c$, all vertices
of $V_i$ have the same colour-$c$ neighbourhood.
Each such neighbourhood consists of three classes $W_{ab}$,
at most one of which receives additional vertices,
so its size is at most
\begin{equation}\label{eq:bipartite-neighbourhoods}
15+\left\lceil\frac{n-60}{4}\right\rceil
=\left\lceil\frac n4\right\rceil.
\end{equation}

Since each $y\in Y$ is joined to the four classes in
distinct colours, the vertices $x,x'\in X$ in any
monochromatic subpath $xyx'$ belong to the same class $V_i$.
Hence every monochromatic path meets at most one $V_i$.
If it meets $V_1$ or $V_2$, all its vertices in $Y$
lie in a single colour neighbourhood of that class.
By \eqref{eq:bipartite-neighbourhoods}, it therefore
contains at most $\lceil n/4\rceil$ vertices of $Y$
and, by alternation, at most $\lceil n/4\rceil+1$
vertices of $X$.

Suppose that at most seven vertex-disjoint monochromatic
paths cover $V(K_{n,n})$. Assign an arbitrary colour to
each singleton path. Since
$|V_1|>2(\lceil n/4\rceil+1)$ and
$|V_2|>\lceil n/4\rceil+1$, at least three paths meet
$V_1$ and at least two meet $V_2$.
Each path meets at most one $V_i$, so covering $V_3$
and $V_4$ requires two further paths.
Thus exactly seven paths are used: three meet $V_1$,
two meet $V_2$, and one meets each of $V_3,V_4$.
The three paths meeting $V_1$ omit a colour $a$,
while the two meeting $V_2$ omit at least two colours,
including some $b\neq a$. None of these five paths
can meet $W_{ab}$, since all edges from $W_{ab}$ to
$V_1,V_2$ have colours $a,b$, respectively.
Each remaining path contains only one vertex of $X$
and hence at most two vertices of $Y$. Together they
cover at most four vertices of $W_{ab}$, contradicting
$|W_{ab}|\geq5$. Therefore at least eight paths are required.
\end{proof}

\section*{Declaration on the Use of Generative AI}
During the preparation of this work, the authors used ChatGPT 6 Pro to assist in discussing proof strategies, checking proofs, and improving exposition.
The authors reviewed and edited the output as needed and take full responsibility for the content of the published article.


\begin{thebibliography}{99}
\bibitem{Allen}
P.~Allen,
\emph{Covering two-edge-coloured complete graphs with two disjoint
monochromatic cycles},
Combin. Probab. Comput. \textbf{17} (2008), no.~4, 471--486.

\bibitem{BQT}
F.~S.~Benevides, A.~L.~Quintino and A.~Talon,
\emph{Partitioning 2-edge-coloured bipartite graphs into monochromatic
cycles},
European J. Combin. \textbf{130} (2025), 104192.

\bibitem{BT}
S.~Bessy and S.~Thomass\'e,
\emph{Partitioning a graph into a cycle and an anticycle, a proof
of Lehel's conjecture},
J. Combin. Theory Ser. B \textbf{100} (2010), no.~2, 176--180.


\bibitem{EGP}
P.~Erd\H{o}s, A.~Gy\'arf\'as and L.~Pyber,
\emph{Vertex coverings by monochromatic cycles and trees},
J. Combin. Theory Ser. B \textbf{51} (1991), no.~1, 90--95.


\bibitem{GG}
L.~Gerencs\'er and A.~Gy\'arf\'as,
\emph{On Ramsey-type problems},
Ann. Univ. Sci. Budapest. E\"otv\"os Sect. Math.
\textbf{10} (1967), 167--170.

\bibitem{G83}
A.~Gy\'arf\'as,
\emph{Vertex coverings by monochromatic paths and cycles},
J. Graph Theory \textbf{7} (1983), no.~1, 131--135.


\bibitem{G}
A.~Gy\'arf\'as,
\emph{Covering complete graphs by monochromatic paths},
in: G.~Hal\'asz and V.~T.~S\'os (eds.),
Irregularities of Partitions,
Algorithms and Combinatorics \textbf{8},
Springer, Berlin, 1989, pp.~89--91.


\bibitem{Gsurvey}
A.~Gy\'arf\'as,
\emph{Vertex covers by monochromatic pieces---a survey of results
and problems},
Discrete Math. \textbf{339} (2016), no.~7, 1970--1977.


\bibitem{GJS}
A.~Gy\'arf\'as, A.~Jagota and R.~H.~Schelp,
\emph{Monochromatic path covers in nearly complete graphs},
J. Combin. Math. Combin. Comput. \textbf{25} (1997), 129--144.

\bibitem{GL}
A.~Gy\'arf\'as and J.~Lehel,
\emph{A Ramsey-type problem in directed and bipartite graphs},
Period. Math. Hungar. \textbf{3} (1973), no.~3--4, 299--304.


\bibitem{GRSS}
A.~Gy\'arf\'as, M.~Ruszink\'o, G.~N.~S\'ark\"ozy and E.~Szemer\'edi,
\emph{An improved bound for the monochromatic cycle partition number},
J. Combin. Theory Ser. B \textbf{96} (2006), no.~6, 855--873.


\bibitem{GRSS3}
A.~Gy\'arf\'as, M.~Ruszink\'o, G.~N.~S\'ark\"ozy and E.~Szemer\'edi,
\emph{Partitioning 3-colored complete graphs into three
monochromatic cycles},
Electron. J. Combin. \textbf{18} (2011), no.~1, Paper P53.

\bibitem{Haxell}
P.~E.~Haxell,
\emph{Partitioning complete bipartite graphs by monochromatic cycles},
J. Combin. Theory Ser. B \textbf{69} (1997), no.~2, 210--218.

\bibitem{LSS}
R.~Lang, O.~Schaudt and M.~Stein,
\emph{Almost partitioning a 3-edge-colored $K_{n,n}$ into five
monochromatic cycles},
SIAM J. Discrete Math. \textbf{31} (2017), no.~2, 1374--1402.

\bibitem{LS}
R.~Lang and M.~Stein,
\emph{Local colourings and monochromatic partitions in complete
bipartite graphs},
European J. Combin. \textbf{60} (2017), 42--54.

\bibitem{Letzter}
S.~Letzter,
\emph{An improvement on \L uczak's connected matchings method},
Bull. Lond. Math. Soc. \textbf{54} (2022), no.~2, 609--623.

\bibitem{LW}
H.~Li and L.~Wang,
\emph{Superlinear lower bounds for monochromatic path partitions}, (2026), arXiv:2608.09895v2.

\bibitem{LRS}
T.~\L uczak, V.~R\"odl and E.~Szemer\'edi,
\emph{Partitioning two-coloured complete graphs into two
monochromatic cycles},
Combin. Probab. Comput. \textbf{7} (1998), no.~4, 423--436.

\bibitem{PRR}
Y.~Peng, V.~R\"odl and A.~Ruci\'nski,
\emph{Holes in graphs},
Electron. J. Combin. \textbf{9} (2002), no.~1, Research Paper 1.

\bibitem{P}
A.~Pokrovskiy,
\emph{Partitioning edge-coloured complete graphs into monochromatic
cycles and paths},
J. Combin. Theory Ser. B \textbf{106} (2014), 70--97.

\bibitem{P23}
A.~Pokrovskiy,
\emph{Partitioning a graph into a cycle and a sparse graph},
Discrete Math. \textbf{346} (2023), no.~1, 113161.

\bibitem{Rado}
R.~Rado,
\emph{Monochromatic paths in graphs},
Ann. Discrete Math. \textbf{3} (1978), 191--194.

\bibitem{Sarkozy}
G.~N.~S\'ark\"ozy,
\emph{Monochromatic cycle partitions of edge-colored graphs},
J. Graph Theory \textbf{66} (2011), no.~1, 57--64.


\bibitem{SS}
O.~Schaudt and M.~Stein,
\emph{Partitioning two-coloured complete multipartite graphs into
monochromatic paths and cycles},
J. Graph Theory \textbf{91} (2019), no.~2, 122--147.


\bibitem{Stein}
M.~Stein,
\emph{Monochromatic paths in 2-edge-coloured graphs and hypergraphs},
Electron. J. Combin. \textbf{30} (2023), no.~1, Paper P1.53.
\end{thebibliography}
\end{document}